\documentclass[reqno]{amsart}
\usepackage{amssymb}
\usepackage[mathscr]{euscript}

\newtheorem{thm}{Theorem}[section]
\newtheorem{cor}[thm]{Corollary}
\newtheorem{lm}[thm]{Lemma}
\newtheorem{con}[thm]{Conjecture}
\newtheorem{prop}[thm]{Proposition}

\newtheorem{exm}[thm]{Example}

\theoremstyle{definition}
\newtheorem{df}[thm]{Definition}

\newcommand{\alg}[1]{{\mathbf #1}}        
\newcommand{\A}{\alg {A}}
\newcommand{\B}{\alg {B}}

\newcommand{\F}{\alg {F}}
\newcommand{\G}{\alg {G}}
\newcommand{\M}{\alg {M}}

\DeclareMathOperator{\V}{V} \DeclareMathOperator{\Id}{Id}
\DeclareMathOperator{\Sub}{\mathbf{CSub}}
\DeclareMathOperator{\Glob}{\mathbf{Cm}}
\DeclareMathOperator{\glob}{\mathcal P}
\def\var{{\mathcal V}}
\def\={ \approx  }

\begin{document}

\author{Kira Adaricheva}
\address{Harold Washington College, 30 East Lake St.,
Chicago, IL 60601, USA} \email{kadaricheva@ccc.edu}

\author{Agata Pilitowska}
\address{Warsaw University of Technology, Plac Politechniki 1,
00-661 Warsaw, Poland} \email{irbis123@o2.pl}

\author{David Stanovsk\'y}
\address{Charles University, Sokolovsk\'a 83, 186 00 Prague, Czech Republic}
\email{stanovsk@karlin.mff.cuni.cz}

\begin{abstract}
Let $\var$ be a variety of algebras. We establish a condition (so
called \emph{generalized entropic property}), equivalent to the
fact that for every algebra $\A\in\var$, the set of all
subalgebras of $\A$ is a subuniverse of the complex algebra of
$\A$. We investigate the relationship between the generalized
entropic property and the entropic law. Further, provided the
generalized entropic property is satisfied in $\var$, we study the
identities satisfied by the complex algebras of subalgebras of
algebras from $\var$.
\end{abstract}

\thanks{\textup{2000} \emph{Mathematics Subject Classification.} 06B20, 06B05}

\thanks{While working on this paper, the authors were supported by the INTAS grant \#03-51-4110.
The third author was also supported by the research project MSM 0021620839 financed by M\v SMT \v CR and by the
grant \#201/05/0002 of the Grant Agency of the Czech Republic.}



\title{On complex algebras of subalgebras}
\maketitle
\begin{center}
\textsl{Dedicated to the 70th birthday of George Gr\"atzer}
\end{center}

\section{Introduction}

For an algebra $\A=(A,F)$, we define \emph{complex operations} on
the set $\glob (A)$ of all non-empty subsets of the set $A$ by
$$f(A_1,\dots,A_n)=\{f(a_1, \dots, a_n): a_i \in A_i\}$$
for every $\emptyset\neq A_1,\dots,A_n\subseteq A$ and every
$n$-ary $f\in F$. The algebra $\Glob \A=(\glob (A),F)$ is called the
\emph{complex algebra} of $\A$.

Complex algebras (called also \emph{globals} or \emph{powers} of
algebras) were studied by several authors, for instance  G.
Gr\"atzer and H. Lakser \cite{GL}, S. Whitney \cite{GW}, A.
Shafaat \cite{Sh}, C. Brink \cite{Br}, I. Bo\v snjak and R.
Madar\'asz \cite {BM}.

The notation of complex operations is used widely. In groups, for
instance, a coset $xN$ is the complex product of the singleton
$\{x\}$ and the subgroup $N$. For a lattice $\mathbf L$, the set
$\Id \mathbf L$ of its ideals forms a lattice under the set
inclusion. If $\mathbf L$ is distributive, then joins and meets in
$\Id \mathbf L$ are precisely the complex operations obtained from
joins and meets of $\mathbf L$, so $\Id \mathbf L$ is a subalgebra
of $\Glob \mathbf L$.

\medskip
Now, consider the set $\Sub \A$ of all (non-empty) subalgebras of
an algebra $\A$. This set may or may not be closed under the
complex operations. For instance, if $\A$ is an abelian group, it
is; however, for most groups, it is not. In the former case, $\Sub
\A$ is a subuniverse of $\Glob \A$ and it will be called the
\emph{complex algebra of subalgebras}. We will say that $\A$
\emph{has the complex algebra of subalgebras} or that $\Sub \A$
\emph{exists}.

Complex algebras of subalgebras were introduced and studied by A.
Romanowska and J.\,D.\,H. Smith in \cite{RS}. A very natural
setting for considering the complex algebras of subalgebras is the
variety of \emph{modes} (idempotent entropic algebras). Research
on complex algebras of submodes was carried out by A. Romanowska
and J.\,D.\,H. Smith in \cite{RS1}, \cite{RS2}, and by the second
author of this paper in \cite{P}, \cite{P1}. In \cite{P2}, the
complex algebras of subalgebras were considered also in the
non-idempotent case.

\medskip
We are concerned with the following question: \emph{In which
varieties does every algebra have the complex algebra of
subalgebras?}

In Section 2 we establish the \emph{generalized entropic property}
for a variety, equivalent to the fact that every algebra has the
complex algebra of subalgebras. The generalized entropic property
appears to be a weak version of the entropic law, so it is natural
to ask about their relationship.

The relationship is investigated in Sections 3 and 4. In general, the generalized entropic property and the
entropic law are not equivalent. We provide
several examples: An idempotent algebra with many binary operations
(Example \ref{rr}), a non-idempotent groupoid (Example \ref{non-idemp1}) or
unary algebras (Example \ref{non-idemp2}). On the other hand,
the generalized entropic property and the entropic law are equivalent under
several additional assumptions, e.g., in groupoids
with a unit element, in commutative idempotent groupoids, or in
idempotent semigroups. We provide several partial results towards the
conjecture that the two conditions are equivalent for idempotent groupoids
(Theorem \ref{idemp_linear} and other).

In Sections 5 and 6, we continue the research started by the
second author in \cite{P1} and investigate which identities are
satisfied by complex algebras of subalgebras. A characterization
of such identities is proved in Theorem \ref{gl for vs}. However,
we are not able to decide the validity of a conjecture stated in
\cite{P1} saying that the variety generated by complex algebras of
subalgebras for algebras from an idempotent variety $\var$
coincides with $\var$ if and only if the latter has a basis of
linear and idempotent identities. We show that a similar statement
for non-idempotent varieties is false, according to Example
\ref{var_example}.

\subsection*{Notation and terminology.}
We denote by $\F_\var(X)$ the free algebra over a set $X$ in a
variety $\var$ and we assume the standard representation of the
free algebra by terms modulo the identities of $\var$. The
notation $t(x_1,\dots,x_n)$ means that the term $t$ contains no
other variables than $x_1,\dots,x_n$ (but not necessarily all of
them) and we say that $t$ is $n$-ary; equivalently, we write
$t\in\F(\{x_1,\dots,x_n\})$.

We call a term $t$ \emph{linear}, if every variable occurs in $t$
at most once. An identity $t\=u$ is called \emph{linear}, if both
terms $t,u$ are linear. An identity $t\=u$ is called
\emph{regular}, if $t,u$ contain the same variables.

An algebra $\A=(A,F)$ is called \emph{entropic} if and only if it satisfies
for every $n$-ary $f\in F$ and $m$-ary $g\in F$ the identity
$$g(f(x_{11},\dots,x_{n1}),\dots,f(x_{1m},\dots,x_{nm}))
\=f(g(x_{11},\dots,x_{1m}),\dots,g(x_{n1},\dots,x_{nm}))$$
(in other words, if all operations of $\A$ commute each other).
Note that a \emph{groupoid}, i.e., a binary algebra, with the
operation denoted usually multiplicatively, is entropic iff it
satisfies the identity
$$xy\cdot uv\=xu\cdot yv,$$
called sometimes the {\it mediality} \cite{JK}. A variety $\var$
is called \emph{entropic} if every algebra in $\var$ is entropic.

An algebra is \emph{idempotent} if each element forms a one-element subalgebra.
Idempotent entropic algebras are called \emph{modes}. The
monograph by A.~Romanowska and J.D.H.~Smith \cite{RS3} provides the most full up-to-date account of results about modes.

\section{generalized entropic property}\label{CompCon}
In this section we introduce and discuss the central notion of
this paper, the generalized entropic property.

\begin{df}\label{generalized entropic property}
We say that a variety $\var$ (respectively, an algebra $\A$)
satisfies the \emph{generalized entropic property} if for every $n$-ary
operation $f$ and $m$-ary operation $g$ of $\var$ (of $\A$), there
exist $m$-ary terms $t_1,\dots,t_n$ such that the identity
$$g(f(x_{11},\dots,x_{n1}),\dots,f(x_{1m},\dots,x_{nm}))
\=f(t_1(x_{11},\dots,x_{1m}),\dots,t_n(x_{n1},\dots,x_{nm}))$$
holds in $\var$ (in $\A$).
\end{df}

For example, a groupoid satisfies the generalized entropic property, if there
are binary terms $t,s$ such that the identity
\begin{gather*}
xy\cdot uv\=t(x,u)s(y,v)\tag{G1}
\end{gather*}
holds. The entropic law is a special case of the generalized entropic property, where the terms
$t_1,\dots,t_n$ are
equal to $g$.

It was proved by T. Evans in \cite{E} that, for a variety $\var$
of groupoids, every groupoid in $\var$ has the complex algebra of
subalgebras iff $\var$ satisfies the identity (G1) for some $s,t$.
We prove the statement for an arbitrary signature. The ``if'' part
of it first appeared in \cite{P1}, where the generalized entropic
property was presented as a ``complex condition''.

\begin{prop}
Every algebra in a variety $\var$ has the complex algebra of
subalgebras i.e., $\Sub\A$ is a subuniverse of $\Glob \A$ for
every $\A\in\var$, if and only if the variety $\var$ satisfies the
generalized entropic property.
\end{prop}

\begin{proof}
First, assume that a variety $\var$ satisfies the generalized entropic property. Let $\A$ be an algebra from the variety $\var$,
$\A_1,\dots,\A_n$ subalgebras of $\A$ and $f$ an $n$-ary operation
of $\A$. We are going to show that $f(A_1,\dots,A_n)$ is closed on
an $m$-ary operation $g$.

Let $x_1,\dots,x_m\in f(A_1,\dots,A_n)$. There exist elements
$a_{ij}\in A_i$, for $1\leq i \leq n$ and $1\leq j \leq m$ such
that $x_j=f(a_{1j},\dots,a_{nj})$. It follows from the generalized entropic property that there exist terms $t_1,\dots,t_n$ such that
\begin{align*}
g(x_1,\dots,x_m)=g(f(a_{11},\dots,a_{n1}),\dots,f(a_{1m},\dots,a_{nm}))=\\
f(t_1(a_{11},\dots,a_{1m}),\dots,t_n(a_{n1},\dots,a_{nm}))\in
f(A_1,\dots,A_n).
\end{align*}
Consequently, $f(A_1,\dots,A_n)$ is a subalgebra of $\A$.

Now, assume that for each algebra $\A\in \var$, the set $\Sub \A$
is closed under complex products. Let $X$ be an infinite set of
variables, let $x_{ij}$, $i=1,\dots,n$, $j=1,\dots,m$, be pairwise
distinct variables from $X$ and let $\F_i$ be the subalgebra of
$\F_\var(X)$ generated by the set $\{x_{ij}\in X\mid
j=1,\dots,m\}$, for every $i=1,\dots,n$. Note that the $\F_i$ are
pairwise disjoint.

For each $n$-ary operation $f$, the set
$f(F_1,\dots,F_n)$ is a subalgebra of $\F_\var(X)$. So for any
$m$-ary operation $g$ and $a_1,\dots,a_m\in f(F_1,\dots,F_n)$, the element
$g(a_1,\dots,a_m)$ belongs to $f(F_1,\dots,F_n)$. Particularly, if
\begin{align*}
a_1&=f(x_{11},\dots,x_{n1}),\\
&\vdots\\
a_m&=f(x_{1m},\dots,x_{nm}),
\end{align*}
then we have
$$g(a_1,\dots,a_m)=g(f(x_{11},\dots,x_{n1}),\dots,f(x_{1m},\dots,x_{nm}))\in f(F_1,\dots,F_n).$$
So there are elements $b_1\in F_1,\dots,b_n\in F_n$ such, that
$g(a_1,\dots,a_m)=f(b_1,\dots,b_n)$. It means that there exist
terms $t_i(x_{i1},\dots,x_{im})$, $i=1,\dots,n$, such that the generalized entropic property
is satisfied in $\F_\var(X)$, and hence in $\var$ too.
\end{proof}

\begin{exm}\label{g1}
{}
\end{exm}
Consider the following 3-element groupoid $\G_1$:
$$\begin{array}{c|ccc} \cdot & a & b & c \\
\hline a & a & c & c\\
b & c & b & c \\
c & a & b & c
\end{array}$$
Notice that $\G_1$ is not entropic, because $c=aa\cdot ba\neq ab\cdot aa=a$.
It is easy to see that $\Sub\G_1$ is a subgroupoid of $\Glob \G_1$, with the multiplication table
$$\begin{array}{c|cccccc} \cdot & \{a\} & \{b\} & \{c\} & \{a,c\} & \{b,c\} & \{a,b,c\} \\
\hline \{a\} & \{a\} & \{c\} & \{c\} & \{a,c\} & \{c\} & \{a,c\} \\
\{b\} & \{c\} & \{b\} & \{c\} & \{c\} & \{b,c\} & \{b,c\} \\
\{c\} & \{a\} & \{b\} & \{c\} & \{a,c\} & \{b,c\} & \{a,b,c\} \\
\{a,c\} & \{a\} & \{b,c\} & \{c\} & \{a,c\} & \{b,c\} & \{a,b,c\} \\
\{b,c\} & \{a,c\} & \{b\} & \{c\} & \{a,c\} & \{b,c\} & \{a,b,c\} \\
\{a,b,c\} & \{a,c\} & \{b,c\} & \{c\} & \{a,c\} & \{b,c\} &
\{a,b,c\}
\end{array}$$
However, there is a groupoid $\F$ in the variety $\V(\G_1)$,
namely $\F=\F_{\V(\G_1)}(x,y,z)$, such that $\Sub \F$ is not a
subgroupoid of $\Glob \F$. To see this, consider the subgroupoid $\A$
of $\F$ generated by $x,y$ and the subgroupoid $\B$ consisting of
$z$. One can check that $A=\{x,y,xy,yx\}$ and that the elements
$((yx)z)y$ and $x\=(xz)x$ are in $(AB)A$. But their product $(((yx)z)y)x$
is not in $(AB)A$. Hence the set $(AB)A$ is not a subgroupoid of
$\F$ and thus $\G_1$ does not satisfy the generalized entropic property.

Later we prove a criterion, Corollary \ref{tr_id}, which shows
that $\G_1$ does not satisfy the generalized entropic property without finding
a particular failure in $\V(\G_1)$.

\medskip
The entropic law is a special case of the generalized entropic property. In the following two sections,
we would like to investigate how far is the generalized entropic property from entropy.

\section{generalized entropic property vs. entropy: the idempotent case}\label{idempotent}


Generally, the generalized entropic property and the entropic law
are not equivalent in idempotent algebras. A counterexample is constructed
in Example \ref{rr}. However, the example has many operations and each
of them is entropic. In the rest of the section, we provide several sufficient conditions
implying that an idempotent groupoid satisfying the generalized entropic property is entropic. The
main result, Theorem \ref{idemp_linear}, is applied several times in the following
propositions and examples.

\begin{exm}\label{rr}
{}
\end{exm}
Let $\mathbf R$ be a ring with a unit, $\G$ a subgroup of the
multiplicative monoid of $\mathbf R$, and $X$ a subset of $G$
closed under conjugation by elements of $X$ and closed under the mapping
$x\mapsto 1-x$, where $-$ is the ring subtraction.

If $\mathbf M$ is a left module over the ring $\mathbf R$, we
define for every element $r\in R$ a binary operation
$\underline{r}:M^2\rightarrow M$ by
$$\underline{r}(x,y)=(1-r)x+ry.$$
Of course, the groupoid $(M,\underline{r})$ is idempotent and
entropic for every $r\in R$. Now, consider the algebra
$\underline{\mathbf M}=(M,\underline{X})$, where
$\underline{X}=\{\underline{r}|r\in X\}$. For every $r,t\in X$, we
put $s_1=(1-r)^{-1}t(1-r)\in X$ and $s_2=r^{-1}tr\in X$ and we get
\begin{gather*}
\underline{t}(\underline{r}(x_1,x_2),
\underline{r}(y_1,y_2))\=(1-t)(1-r)x_1+(1-t)rx_2+t(1-r)y_1+try_2\=\\
\=(1-r)(1-s_1)x_1+r(1-s_2)x_2+(1-r)s_1y_1+rs_2y_2\=
\underline{r}(\underline{s_1}(x_1,y_1),\underline{s_2}(x_2,y_2)).
\end{gather*}
So the algebra $\underline{\mathbf M}$ satisfies the generalized entropic property. On the other hand, it is entropic, iff $rt=tr$ for all
$r,t\in X$. To check this put $x_1=y_1=y_2=0$ and $x_2=1$ in the
previous identity.

For example, if $\mathbf R$ is a non-commutative division ring (a
skew field), $\G$ is its multiplicative group and
$X=R\smallsetminus\{0,1\}$, then $\underline{\mathbf M}$ is a
non-entropic idempotent algebra satisfying the generalized entropic property.

However, it is infinite, with infinitely many (binary) operations.
To get a finite example, we need a more elaborate setting.

Let $\mathbf R$ be the ring of all $2\times 2$ matrices over a
field $\mathbb F$, $\G$ the subgroup of all matrices with
determinant 1 and $X$ the subset of all matrices with trace 1. It
is well known that traces are invariant under conjugation and it
is easy to check that $X$ is closed under the mapping
$x\mapsto 1-x$. Let $\M$ be a two-dimensional vector space over
$\mathbb F$, considered as a module over $\mathbf R$.

If $\mathbb F=\textup{GF}(2)$ then $X$ has only two elements and they
commute. If $\mathbb F=\textup{GF}(3)$, then $X$ has nine elements and some
of them do not commute, so we get a 9-element non-entropic
idempotent algebra $\M_9=(M_9,\underline{X})$ with 9 binary
operations satisfying the generalized entropic property. In fact, the algebra
$(M_9,\underline{X'})$, where
$X'=X\smallsetminus\{\left(\smallmatrix 2&0 \\ 0 &
2\endsmallmatrix\right)\}$, has the same properties.

Finally, we note that similar examples can be obtained with operations of an arbitrary
arity $n\geq2$; consider the operations
$$\underline{(r_2,\dots,r_n)}(x_1,\dots,x_n)=(1-r_2-\dots-r_n)x_1+r_2x_2+\dots+r_nx_n.$$

Notice that the algebra $(M,\underline r)$ is entropic, for any $r$. So, one might think about the following conjecture:

\begin{con}\label{conjidemp}
Every idempotent algebra $(A,f)$ with the generalized entropic property is entropic.
\end{con}

In the sequel, we prove several special cases of the conjecture for groupoids.

As we previously noticed, the generalized entropic property in groupoids is
equivalent to the following statement: There are binary terms $t$
and $s$ such that the identity
\begin{gather*}
xy\cdot uv\=t(x,u)s(y,v)\tag{G1}
\end{gather*}
holds. An immediate consequence of the generalized entropic property in
idempotent groupoids are the following important identities that
can be treated as the laws of \emph{pseudo-distributivity}:
\begin{gather*}
xy\cdot xz\=xs(y, z)\tag{G2},\\
yx\cdot zx\=t(y, z)x\tag{G3},\\
x\cdot yz\=t(x, y)s(x, z)\tag{G4},\\
yz\cdot x\=t(y, x)s(z, x)\tag{G5}.
\end{gather*}
(G2) states that, for every $a$, the \emph{left translation} $L_a:x\mapsto ax$ is a homomorphism $(G,s)\to(G,\cdot)$
and (G3) states that the \emph{right translation} $R_a:x\mapsto xa$ is a homomorphism $(G,t)\to(G,\cdot)$.

The main partial result towards Conjecture \ref{conjidemp} is the following theorem.

\begin{thm}\label{idemp_linear}
If an idempotent groupoid $\G$ satisfies the generalized entropic property for
some terms $t,s$ and at least one of $t,s$ is linear, then $\G$ is
entropic.
\end{thm}

\begin{proof}
If $t$ is linear, one of Lemmas \ref{lem1}--\ref{lem4}, applies.
If $s$ is linear, consider the dual groupoid $\G^\partial$
(with the operation defined by $x\bullet y=yx$);
this groupoid satisfies the generalized entropic property with the role of
$t,s$ interchanged, hence both $\G^\partial$ and $\G$ are entropic
by one of Lemmas \ref{lem1}--\ref{lem4}; note that entropy is a self-dual identity.
\end{proof}

\begin{lm}\label{lem1}
If an idempotent groupoid $\G$ satisfies the generalized entropic property for the term $t(x,y)=x$
and an arbitrary term $s$, then $\G$ is entropic.
\end{lm}

\begin{proof}
The generalized entropic property states that $xy\cdot uv\=xs(y,v)$. Since the
value of $xy\cdot uv$ does not depend on $u$, we have $xy\cdot
uv\=xy\cdot vv\=xy\cdot v$. Hence, with $x=y$, we obtain $x\cdot
uv\=xv$. Applying this identity to the term $xs(y,v)$, we get
$xs(y,v)\=xw$, where $w\in\{y,v\}$ is the rightmost variable in
the term $s(y,v)$. So, we have $xy\cdot uv\=xy\cdot
v\=xs(y,v)\=xw$. If $w=y$, then $xv\=xx\cdot vv\=xx\=x$ by identifying:
$x=y$ and $u=v$. Thus the entropy holds. If $w=v$ then $xy\cdot
uv$ does not depend on $y$ and $u$, hence we can interchange them
and the entropy holds again.
\end{proof}

\begin{lm}
If an idempotent groupoid $\G$ satisfies the generalized entropic property for the term $t(x,y)=y$
and an arbitrary term $s$, then $\G$ is entropic.
\end{lm}

\begin{proof}
The generalized entropic property says that $xy\cdot uv\=us(y,v)$. Since the
value of $xy\cdot uv$ does not depend on $x$, we have $xy\cdot
uv\=yy\cdot uv\=y\cdot uv$. Hence, with $u=v$ we obtain $xy\cdot
u\=yu$. Applying this identity to the term $s(x,y)z$, we get
$s(x,y)z\=wz$, where $w\in\{x,y\}$ is the rightmost variable in
the term $s(x,y)$. So, we have
$s(x,y)\=s(x,y)s(x,y)\=ws(x,y)=t(x,w)s(x,y)$, and thus $s(x,y)\=x\cdot wy$ by
the generalized entropic property. So we may assume that the rightmost
variable of $s$ is $y$, i.e., $w=y$. Consequently, $s(x,y)\=xy$ and thus
$xy\cdot uv\=u\cdot yv\=y\cdot uv\=xu\cdot yv$ by (G1), (G4) and
(G1).
\end{proof}

\begin{lm}
If an idempotent groupoid $\G$ satisfies the generalized entropic property for the term $t(x,y)=xy$
and an arbitrary term $s$, then $\G$ is entropic.
\end{lm}

\begin{proof}
Note that (G3) is the right distributivity. Hence
$r(x,y)z\=r(xz,yz)$, for every term $r$.

\noindent
{\sl Claim 1.} $s(x,z)\cdot xz\=s(x,z)$.

Using right distributivity in $s$, twice the generalized entropic property and
again the right distributivity in $s$, we obtain
\begin{align*}
s(x,z)\cdot
xz&\=s(x\cdot xz,z\cdot xz)\=s(xs(x,z),zx\cdot z))\=s(xs(x,z),zs(x,z))\\
&\=s(x,z)s(x,z)\=s(x,z).
\end{align*}

\noindent
{\sl Claim 2.} $xs(y,z)\=x\cdot yz$.

Using several times the generalized entropic property and the idempotent law,
we get
\begin{align*}
xs(y,z)&\=xy\cdot xz\=(xy)(xz\cdot xz)\=(x\cdot
xz)s(y,xz)\=(xs(x,z))s(y,xz)\\
&\=(xy)(s(x,z)\cdot
xz)\=(xy)s(x,z)\=x\cdot yz,
\end{align*}
where the last but one equality follows from
Claim 1.

Finally, it follows from Claim 2 that $xy\cdot uv\=xu\cdot s(y,v)\=xu\cdot yv$.
\end{proof}

\begin{lm}\label{lem4}
If an idempotent groupoid $\G$ satisfies the generalized entropic property for the term $t(x,y)=yx$
and an arbitrary term $s$, then $\G$ is entropic.
\end{lm}

\begin{proof}
Note that (G3) is read as $xy\cdot z\=yz\cdot xz$ that can be
treated as the {\it right anti-distributivity}. One can check by
induction that $r(x,y)z\=r^\partial(xz,yz)$ for every term $r$,
where $r^\partial$ denotes the term dual to $r$ (this is the term
that results when reading $r$ from right to left; inductively,
$x^\partial=x$ and $(r_1r_2)^\partial=r_2^\partial r_1^\partial$).

\noindent
{\sl Claim 1.} $s(x,z)\cdot xz\=s(x,z)$.

Using the right anti-distributivity in $s$, then three times the
generalized entropic property and again the right anti-distributivity in $s$, we
get
\begin{align*}
s(x,z)\cdot xz&\=s^\partial(x\cdot xz,z\cdot
xz)\=s^\partial(xs(x,z),xz\cdot z))\=s^\partial(xs(x,z),zx\cdot z))\\
&\=s^\partial(xs(x,z),zs(x,z))\=s(x,z)s(x,z)\=s(x,z).
\end{align*}

\noindent
{\sl Claim 2.} $s(x,y)\=xy$.

Using twice Claim 1 and three times the generalized entropic property, we obtain
\begin{align*}
s(x,y)&\=s(x,y)(xy)\=(s(x,y)\cdot
xy)(xy)\=(xs(x,y))s(xy,y)\=(x\cdot xy)s(xy,y)\\
&\=(xy\cdot xy)(xy)\=xy.
\end{align*}

Hence the groupoid $\G$ satisfies $xy\cdot uv\=ux\cdot yv$.
Consider the dual groupoid $\G^\partial$;
it satisfies $xy\cdot uv\=xu\cdot vy$ and thus it is entropic by the preceding lemma.
Since entropy is a self-dual identity, $\G$ is entropic too.
\end{proof}

Theorem \ref{idemp_linear} has several interesting consequences.

\begin{cor}\label{free2}
Let $\var$ be a variety of idempotent groupoids such that every binary term is equivalent
to a linear term in $\var$. If $\var$ satisfies the generalized entropic property, then $\var$ is entropic.
\end{cor}

All groupoids with the property that every binary term is
equivalent to a linear term were characterized by J. Dudek
\cite{D}, see also \cite{DJMMS}. The groupoid $\G_1$ from Example
\ref{g1} can be found in the list of these groupoids.

\begin{cor}\label{tr_id}
Let $\G$ be an idempotent groupoid with a one-sided unit, i.e., $e\in G$ such
that $ex=x$ for all $x\in G$, or $xe=x$ for all $x\in G$.
If $\G$ satisfies the generalized entropic property, then it is entropic.
\end{cor}

\begin{proof}
Assume that $e$ is a left unit in $\G$. Then
$$xy\=ex\cdot ey\=t(e,e)s(x,y)\=es(x,y)\=s(x,y)$$
in $\G$ and thus Theorem \ref{idemp_linear} applies. If $e$ is a right unit,
proceed dually.
\end{proof}

For example, the element $c$ is a left unit in the groupoid $\G_1$ from Example \ref{g1}.
Since $\G_1$ is non-entropic, it cannot satisfy the generalized entropic property.

The following observation will also become useful in the sequel.

\begin{lm}\label{same}
If an idempotent algebra $\A=(A,F)$ satisfies the generalized entropic property such
that, for each pair $f,g\in F$, the terms $t_1,\dots,t_n$ are equal,
then $\A$ is entropic.
\end{lm}

\begin{proof}
Let $t=t_1=\dots=t_n$. Then
\begin{align*}
g(x_1,\dots,x_m)&\=g(f(x_{1},\dots,x_{1}),\dots,f(x_{m},\dots,x_{m}))\\
&\=f(t(x_{1},\dots,x_{m}),\dots,t(x_{1},\dots,x_{m}))\=t(x_{1},\dots,x_{m}).
\end{align*}
\end{proof}

Now, we apply our previous results in several well-known classes
of groupoids.

\begin{prop}\label{bands}
An idempotent semigroup \textup{(}i.e., a band\textup{)} satisfying the generalized entropic property is entropic.
\end{prop}

\begin{proof}
In bands, any binary term is equivalent to one of $x$, $y$, $xy$,
$yx$, $xyx$, $yxy$: by the idempotency, neither a variable can
appear at two consecutive places, nor $xy$ can appear more then
once in a row. So, if $t$ or $s$ is equivalent to one of the first
four (linear) terms, we can apply Theorem \ref{idemp_linear}. If
$t,s$ are equivalent to the same term, then we can use Lemma
\ref{same}. Hence we are left with two cases:
$$xyuv\=uxuyvy \qquad\text{and}\qquad xyuv\=xuxvyv.$$

If the first identity holds, then we get $xv\=xvx$ by substitution $x=y=u$,
and $xv\=vxv$ by substitution $y=u=v$. So, we have the
commutativity, hence the entropy follows.

In the latter case, $xuw\=xxuw\=xuxwxw\=xuxw$, where the last
equality follows from the idempotency, and, similarly,
$wvyv\=wyv$. Thus $xuxvyv$ is equal to
$(xux)(vyv)\=(xu)(vyv)\=(xu)(yv)$.
\end{proof}

\begin{prop}
An idempotent commutative groupoid satisfying the generalized entropic property is entropic.
\end{prop}

\begin{proof}
Using (G3), the commutativity and (G2), we obtain
$$t(x,y)z\=xz\cdot yz\=zx\cdot zy\=zs(x,y)\=s(x,y)z.$$
Consequently, $$s(x,u) \= s(x,u)s(x,u) \= t(x,u)s(x,u) \= xx\cdot
uu \= xu.$$ Similarly for $t$.
\end{proof}

A groupoid $\G$ is called \emph{left} (respectively, \emph{right})
\emph{cancellative}, if $zx=zy$ implies $x=y$ ($xz=yz$ implies
$x=y$), for all $x,y,z\in G$. For instance, quasigroups are both
left and right cancellative.

\begin{prop}
An idempotent left or right cancellative groupoid satisfying
the generalized entropic property is entropic.
\end{prop}

\begin{proof}
Assume the left cancellativity. Then $x\cdot xy\=xx\cdot
xy\=t(x,x)s(x,y)\=x\cdot s(x,y)$ and so by the left cancellativity
we get $s(x,y)\=xy$. Now apply Theorem~\ref{idemp_linear}. In the
case of the right cancellativity proceed dually.
\end{proof}

Next, we apply Corollary \ref{free2} to show that the generalized entropic property
fails in the varieties generated by all graph algebras and by all equivalence algebras, although
every graph algebra and every equivalence algebra has the complex algebra of subalgebras.

\begin{exm}\label{equiv-alg}
{}
\end{exm}
Let $A$ be a set and let $\alpha \subseteq A\times A$ be an
equivalence relation on $A$. The \emph{equivalence algebra}
$\A(\alpha)$ is a groupoid with the multiplication defined as
follows (see, for example, \cite{JM}): $$x\cdot
y=\left\{\begin{array}{ll}x, & {\rm if} \; \; (x,y)\in \alpha,
\\ y, & {\rm otherwise.}\end{array}\right.$$
It is easy to see that a homomorphic image and a subalgebra of an
equivalence algebra is again an equivalence algebra. In fact, any
subset of an equivalence algebra is a subalgebra. Hence, every
equivalence algebra has the complex algebra of subalgebras.

However, consider the variety $\mathcal E$ generated by
equivalence algebras. It is not entropic, since in the equivalence
algebra on the set $\{a,b,c\}$, corresponding to the equivalence
with two blocks $\{a,b\}$ and $\{c\}$, we have
$a=(ca)b\neq(cb)(ab)=b.$ It is not difficult to check that the
two-generated free algebra in $\mathcal E$ has only four elements:
$x,y,xy,yx$. Hence, by Corollary \ref{free2}, the variety
$\mathcal E$ does not satisfy the generalized entropic property.

\begin{exm}\label{4-alg}
{}
\end{exm}
Let $G=(V,E)$ be a graph with a set $V$ of vertices and a set
$E\subseteq V\times V$ of edges. Its \emph{graph algebra}
$\A(G)=(V\cup \{0\},\cdot)$ is a groupoid with the multiplication
defined as follows:
$$x\cdot
y=\left\{\begin{array}{ll}x, & {\rm if} \; \; (x,y)\in E,
\\ 0, & {\rm otherwise.}\end{array}\right.$$
As shown in \cite{O}, any subalgebra with $0$ and any homomorphic
image of a graph algebra is a graph algebra.

In fact, any subset with $0$ of a graph algebra is clearly a
subalgebra. Moreover, for any two subalgebras $\A$ and $\B$ of the
graph algebra $\A(G)$, if all $a\in \A$ and $b\in \B$ are connected
by the edge, then $\A\B=\A$. On the other side, if there are such $a\in
\A$ and $b\in \B$ that $(a,b)$ is not in $E$, then $0\in \A\B$.
Thus every graph algebra has the complex algebra of subalgebras.

However, consider the variety $\mathcal G_I$ generated by
idempotent graph algebras. It is not entropic, since in the graph
algebra corresponding to the graph

\noindent
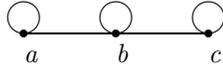
\centerline{
\begin {picture}(150,40)(0,0)
\put(55,35) {} \put(55,15){\circle*{3}}
\put(55,21){\circle{12}} \put(55,4){{\it a}}
\put(55,15){\line(1,0){35}} \put(90,21){\circle{12}}
\put(90,15){\circle*{3}} \put(90,4){{\it b}}
\put(90,15){\line(1,0){35}} \put(125,15){\circle*{3}}
\put(125,21){\circle{12}} \put(125,4){{\it c}}
\end{picture}}
we have $b=(bc)a\neq (ba)(ca)=0$. Similarly to Example
\ref{equiv-alg}, the two-generated free algebra in $\mathcal G_I$
has only four elements: $x,y,xy,yx$. By Corollary \ref{free2}, the
variety $\mathcal G_I$ does not satisfy the generalized entropic property.

\medskip
We finish this section with an observation.

\begin{prop}\label{worthless}
Every idempotent groupoid with the generalized entropic property satisfies the identity
$$xy\cdot uv\=(xy\cdot uy)(xv\cdot uv)\=(xy\cdot xv)(uy\cdot uv).$$
\end{prop}

\begin{proof}
Using (G1) and (G2) we obtain $xy\cdot uv\=t(x,u)s(y,v)\=t(x,u)y\cdot t(x,u)v$ and
now the first identity follows from (G3).
Similarly, using (G1), (G3) and (G2) we obtain
$xy\cdot uv\=t(x,u)s(y,v)\=xs(y,v)\cdot us(y,v)\=(xy\cdot xv)(uy\cdot uv)$.
\end{proof}

The converse is false. It can be checked that the following groupoid $\G_2$
$$\begin{array}{c|ccc} \cdot & a & b & c \\
\hline a & a & b & a\\
b & c & b & c \\
c & c & b & c
\end{array}$$
satisfies the identities from Proposition \ref{worthless},
but it fails the generalized entropic property.

\section{generalized entropic property vs. entropy: the non-idempotent case}\label{non-idempotent}

We start with an observation that the generalized entropic property
and the entropic law are generally inequivalent for non-idempotent
groupoids.

\begin{exm}\label{non-idemp1}
{}
\end{exm}
Let $\var_A$ denote the variety of 
groupoids satisfying the identity $$(x_1x_2)(x_3x_4)\=(x_3x_1)(x_2x_4).$$
Clearly, the generalized entropic property holds in $\var_A$.
However, it follows from Lemma \ref{l_reverse} that $\var_A$ is not entropic:
in our case $A=\{(1,3,2)\}$, and the
subgroup generated by $A$ in the symmetric group $\alg S_4$ does not contain the transposition $(2,3)$.

\begin{lm}\label{l_reverse}
Let $A\subseteq S_4$ be a set of permutations on four elements and let $\var_A$ be
the variety of groupoids satisfying the identities
$$x_1x_2\cdot x_3x_4\=x_{\pi1}x_{\pi2}\cdot x_{\pi3}x_{\pi4}$$
for every $\pi\in A$. Then $\var_A$ is entropic, iff the transposition $(2,3)$ is
in the subgroup generated by $A$ in $\alg S_4$.
\end{lm}

\begin{proof}
Generally, two terms $p,q$ are equivalent in a variety $\var$, iff
there is a sequence $p=w_1,w_2,\dots$, $w_n=q$ such that, for each
$i$, $w_i$ has a subterm which is a substitution instance of some
term which appears in an equation $\varepsilon$ from the base of
$\var$, and $w_{i+1}$ is derived from $w_i$ by replacing this
subterm by the same substitution instance of the other side of
$\varepsilon$. In $\var_A$, starting with the term $x_1x_2\cdot
x_3x_4$, we cannot make any proper substitution, hence $w_{i+1}$
is always obtained from $w_i$ by permuting variables of $w_i$ by
some $\pi\in A$. Hence, if we use permutations
$\pi_1,\dots,\pi_{n-1}$, we arrive in the term
$x_{\pi1}x_{\pi2}\cdot x_{\pi3}x_{\pi4}$, where
$\pi=\pi_{n-1}\cdots\pi_1$. So, the entropy can be obtained iff
$(2,3)$ is generated by permutations from $A$.
\end{proof}

\medskip
The two conditions are inequivalent also for unary algebras.
(They haven't appeared in the previous section, because idempotency is a rather trivial property in there.)

\begin{exm}\label{non-idemp2}
{}
\end{exm}
Let $\A=(A,F)$ be a unary algebra, i.e., $F$ contains only unary
operations. Clearly, $\A$ is entropic iff $fg\=gf$, for all $f,g\in
F$.

Let $B$ be a subset of the symmetric group over a set $X$ such
that $B=B^{-1}$. Put $\B=(X,\{f:f\in B\})$. Then for every $f,g\in
B$ we can always find a term $t$ such that $fg = gt$ (namely, $t =
g^{-1}fg$), so $\B$ satisfies the generalized entropic property. On the other
hand, if $fg\neq gf$ for at least one pair $f,g\in B$, then $\B$
is not entropic.

\bigskip
On the other hand, there are several important classes, where the
generalized entropic property is equivalent with the entropic law, regardless
idempotency. For instance, this is true for groupoids with a unit
element. The following statement covers a more general setting. We
say that an element $e$ is a \emph{unit} for an operation
$f$, if
\begin{gather*}
f(x,e,\dots,e)\=f(e,x,e,\dots,e)\=\dots\=f(e,\dots,e,x)\=x
\end{gather*}
for every $x\in A$. We say that $e$ is a \emph{unit} for an
algebra $(A,F)$, if it is a unit for each operation $f\in F$.

\begin{lm}\label{unitlm}
Let $\A=(A,F)$ be an algebra with a one-element subalgebra $\{e\}$
and assume that $e$ is a unit for an $n$-ary operation $f\in F$.
If $(A,F)$ satisfies the generalized entropic property, then $f$ commutes
with each operation $g\in F$.
\end{lm}

\begin{proof}
The generalized entropic property says that
$$g(f(x_{11},\dots,x_{n1}),\dots,f(x_{1m},\dots,x_{nm}))\=f(t_1(x_{11},\dots,x_{1m}),
\dots,t_n(x_{n1},\dots,x_{nm}))$$ for some terms $t_1,\dots,t_n$.
Hence
\begin{align*}
g(x_1,\dots,x_m)&\=g(f(e,\dots,x_1,\dots,e),f(e,\dots,x_2,\dots,e),
\dots,f(e,\dots,x_m,\dots,e)) \\
&\=f(t_1(e,\dots,e),\dots,t_k(x_1,\dots,x_m),\dots,t_n(e,\dots,e))\\
&\=f(e,\dots,t_k(x_1,\dots,x_m),\dots,e)\=t_k(x_1,\dots,x_m),
\end{align*}
 for every $k \leq n$. So
 $$g(f(x_{11},\dots,x_{n1}),\dots,f(x_{1m},\dots,x_{nm}))
 \=f(g(x_{11},\dots,x_{1m})\dots,g(x_{n1},\dots,x_{nm})).$$
\end{proof}

As a consequence, we get

\begin{prop}\label{unit}
Let $\A$ be an algebra with a unit element $e$.
If $\A$ satisfies the generalized entropic property, then $\A$ is entropic.
\end{prop}

Adjoining an outside unit element is quite a standard operation when dealing with algebras. The following example shows that such an extended algebra may fail the generalized entropic property.
\begin{exm}
{}
\end{exm}
Consider a groupoid $\G$ satisfying the generalized entropic property and
possessing elements $a,b$ such that $ab\neq ba$. Let $\G^\ast$
denote the groupoid obtained from $\G$ by adjoining a unit element
$e$. Then $\G^\ast$ is not entropic, because $ea\cdot be=ab\neq
ba=eb\cdot ae.$ Hence, although $\G$ itself satisfies the generalized entropic property, by Proposition \ref{unit} the groupoid $\G^\ast$ does
not.

\medskip
A \emph{loop} is an algebra $\A=(A,\cdot,/,\backslash,e)$ such that
the identities
\begin{gather*}
x\backslash(xy)\=y,\qquad(yx)/x\= y,\\
x(x\backslash y)\=y,\qquad(y/x)x\= y,\\
xe\=ex\=x
\end{gather*}
hold in $\A$. In other words, loops can be considered as
``non-associative groups''. On the other hand, groups can be
regarded as loops with $x/y=xy^{-1}$ and $y\backslash x=y^{-1}x$.

\begin{prop}\label{loops}
Let $\var$ be a variety of loops. The following conditions are
equivalent:
\begin{enumerate}
\item $\var$ satisfies the generalized entropic property;
\item $\var$ is entropic;
\item $\var$ is a variety of abelian groups.
\end{enumerate}
\end{prop}

\begin {proof}
(1) $\Rightarrow$ (2). Let $\A\in\var$. It follows from Lemma \ref{unitlm} that $(A,\cdot)$ is entropic.
And it is easy to check that, for any loop $\A$, if $(A,\cdot)$ is entropic, then $\A$ is entropic.

(2) $\Rightarrow$ (3).
If $\A$ is an entropic loop and $x,y,z\in A$, then
$$xy\cdot z=xy\cdot ez=xe\cdot yz=x\cdot yz$$ (hence $\A$ is a group) and
$$xy=(xy)(x^{-1}x)=(xx^{-1})(yx)=yx.$$

(3) $\Rightarrow$ (1). It is well known that the product of two subgroups is a subgroup.
\end{proof}

\medskip
We finish this section with a result on commutative groupoids.
A term $r$ is called \emph{$\G$-symmetric}, if $\G$ satisfies $r(x,y)\=r(y,x)$.

\begin{prop}
If a commutative groupoid $\G$ satisfies the generalized entropic property for some terms $t,s$ and
at least one of $t,s$ is linear or $\G$-symmetric, then $\G$ is entropic.
\end{prop}

\begin{proof}
Because of commutativity, we can assume that the linear or
$\G$-symmetric term is $t$. If $t$ is $\G$-symmetric, then, using several
times the commutativity and the generalized entropic property, we get
$$xy\cdot uv\=yx\cdot uv\=t(y,u)\cdot s(x,v)\=t(u,y)\cdot s(x,v)\=ux\cdot yv\=xu\cdot yv.$$
If $t$ is linear, then either $t(x,y)\in\{xy,yx\}$ (so $t$ is
$\G$-symmetric and the first case applies), or $t(x,y)=x$, or $t(x,y)=y$. First, assume $xy\cdot
uv\=xs(y,v)$. Consequently, the term $xy\cdot uv$ does not depend
on $u$ and we can compute using the commutativity:
$$xy\cdot uv\=xy\cdot yv\=yv\cdot xy\=yv\cdot yx\=yv\cdot ux\=ux\cdot yv\=xu\cdot yv.$$
Next, if $xy\cdot uv\=us(y,v)$, then $xy\cdot uv$ does not depend on $x$ and a similar computation does the job.
\end{proof}

\section{Identities in complex algebras of subalgebras}

Let $\var$ be a variety. We will denote by $\Glob\var$ the variety
generated by complex algebras of algebras in $\var$, i.e.,
$$\Glob\var=\V(\{\Glob \A:\A\in\var\}).$$
Further, if $\var$ satisfies the generalized entropic property, we let
$\Sub\var$ be the variety generated by complex algebras of
subalgebras of algebras in $\var$, i.e.,
$$\Sub\var=\V(\{\Sub \A :\A\in\var\}).$$
Evidently, $\Sub\var\subseteq\Glob\var$, because $\Sub \A$ is a subalgebra
of $\Glob \A$. Also $\var\subseteq\Glob\var$, because every algebra $\A$ can
be embedded into $\Glob \A$ by $x\mapsto\{x\}$. And if $\var$ is idempotent,
then $\var\subseteq\Sub\var$, by the same embedding. On the other hand,
we do not have $\var\subseteq\Sub\var$ in general, for instance, for the
variety of abelian groups ($\Sub\var$ is defined due to
Proposition \ref{loops}), because in this case $\Sub\var$ is idempotent,
while $\var$ is not.

In \cite{GL}, G.~Gr\"atzer and H.~Lakser proved the following theorem.

\begin{thm}\label{gl}
Let $\var$ be a variety. Then $\Glob\var$ satisfies precisely those identities
resulting through identification of variables from the linear
identities true in $\var$.
\end{thm}

\begin{cor}\label{glcor}
Let $\var$ be a variety. Then
$\var=\Glob\var$, if and only if $\var$ has a base consisting of linear identities.
\end{cor}

We investigate the question raised in \cite{P}:
{\it What are the identities satisfied by $\Sub\var$} (provided it
is defined)? In particular, {\it when $\var=\Sub\var$?}

It follows from Theorem \ref{gl} that $\Sub\var$ satisfies the linear
identities valid in $\var$. However, if $\var$ is idempotent, $\Sub\var$ is
also idempotent, and the idempotency is not linear.
Moreover $\Sub\var$ can still be idempotent, while $\var$
is not---recall the example with abelian groups.
We are going to prove an analogue of Theorem \ref{gl}, characterizing
the identities satisfied by $\Sub\var$. First, we have to introduce the notion
of a \emph{semilinear precursor}.

An identity $t\=s$ is called \emph{semilinear}, if at least one of
the terms $t,s$ is linear. The \emph{linearization} of a term
$t(x_1,\dots,x_n)$ is the term $t^\ast$, resulting from $t$ by replacement
of the $j$-th occurence of a variable $x_i$ by the variable $x_{ij}$, for all $1\leq i\leq n$ and
$1\leq j\leq k_i$, where $k_i$ is the number of occurences of the variable $x_i$ in $t$.

Let $t,s$ be terms and let $k_i, l_i$ denote the number of occurences of the variable $x_i$ in $t, s$.
If $x_i$ does not occur in the term $t$, we redefine $k_i=1$.

The identity $t^\ast\=\tilde s$ is called a \emph{semilinear precursor} for the (ordered) pair $(t,s)$,
if there are terms $r_{ij}(x_{i1},\dots,x_{ik_i})$, $1\leq i\leq n$, $1\leq j\leq l_i$ such that
$$\tilde s=s^\ast(r_{11}(\overline{x_1}),\dots,r_{1l_1}(\overline{x_1}),\dots,r_{n1}(\overline{x_n}),\dots,r_{nl_n}(\overline{x_n}))$$
(where $\overline{x_i}$ denotes the tuple $(x_{i1},\dots,x_{ik_i})$).
For example, the semilinear precursors for the pair $(xy\cdot xz,yz\cdot x)$ are precisely the identities
of the form $x_1y\cdot x_2z\=p(y)q(z)\cdot r(x_1,x_2)$, where $p,q$ are unary terms and $r$ is a binary term.
The semilinear precursors for the pair $(yz\cdot x,xy\cdot xz)$ are precisely the identities
of the form $yz\cdot x\=p_1(x)q(y)\cdot p_2(x)r(z)$, where $p_1,p_2,q,r$ are unary terms. (In both examples, instead
of double indices we used different letters for variables and terms.)

Indeed, the identity $t\=s$ results from any of its semilinear precursors through identification of the variables
$x_{i1},\dots,x_{ik_i}$ and replacement of the unary subterms $r_{ij}(x_i,\dots,x_i)$ by a single variable.
In particular, the identity $t\=s$ is a consequence of each semilinear precursor for the pair $(t,s)$ and idempotency.

\begin{thm}\label{gl for vs}
Let $\var$ be a variety satisfying the generalized entropic property. Then
$\Sub\var$ satisfies the identity $t\=s$, if and only if there are
semilinear precursors for the pair $(t,s)$ and for the pair
$(s,t)$, both satisfied in $\var$.
\end{thm}

\begin{proof}
First, assume that $t(x_1,\dots,x_n)\=s(x_1,\dots,x_n)$ holds in $\Sub\var$ and
denote $k_i,l_i$ the number of occurences of the variable $x_i$ in $t,s$.
Again, if $x_i$ does not occur in the term $t$, we redefine $k_i=1$.

Let $\A_i$, $i=1,\dots,n$, be the subalgebra generated by the set
$\{x_{i1},\dots,x_{ik_i}\}$ in $\F_\var(X)$, the free algebra in $\var$
over the set $X=\{x_{ij}:1\leq i\leq n,1\leq j\leq k_i\}$.
Since $t(A_1,\dots,A_n)=s(A_1,\dots,A_n)$, we have
$$t^\ast(x_{11},\dots,x_{1k_1},\dots,x_{n1},\dots,x_{nk_n})\in s(A_1,\dots, A_n).$$
It means that there are terms $r_{ij}(x_{i1},\dots,x_{ik_i})\in A_i$, $1\leq i\leq n$, $1\leq j\leq l_i$ such that
\begin{gather*}
t^\ast(x_{11},\dots,x_{1k_1},\dots,x_{n1},\dots,x_{nk_n})\=\\
s^\ast(r_{11}(\overline{x_1}),\dots,r_{1l_1}(\overline{x_1}),\dots,r_{n1}(\overline{x_n}),\dots,r_{nl_n}(\overline{x_n})).
\end{gather*}
In other words, the above identity is a semilinear precursor for the pair $(t,s)$ and
it is satisfied in $\var$, because the identity holds in a free algebra.
To get a semilinear precursor for the pair $(s,t)$, consider the same procedure with the role of $t,s$ interchanged.

Now we prove the converse. Let $t(x_1,\dots,x_n),s(x_1,\dots,x_n)$ be terms and assume there are semilinear precursors
$t^\ast\=\tilde{s}$ and $s^\ast\=\tilde{t}$ satisfied in $\var$. Let $\A\in\var$
and take arbitrary subalgebras $\A_1,\dots,\A_n$ of $\A$. To prove the inclusion
$t(A_1,\dots,A_n)\subseteq s(A_1,\dots,A_n)$, let $a\in t(A_1,\dots,A_n)$. It means that there are
$a_{i1},\dots,a_{ik_i}\in A_i$ ($i=1,\dots,n$) such that
$$a=t^\ast(a_{11},\dots,a_{nk_n}).$$
The algebra $\A$ satisfies $t^\ast\=\tilde s$, so
$$a=s^\ast(r_{11}(\overline{a_1}),\dots,r_{1l_1}(\overline{a_1}),\dots,r_{n1}(\overline{a_n}),\dots,r_{nl_n}(\overline{a_n})).$$
Since $r_{ij}(a_{i1},\dots,a_{ik_i})\in A_i$ for every $i,j$, we see that $$a\in s^\ast(A_1,\dots,A_1,\dots,A_n\dots,A_n)=s(A_1,\dots,A_n).$$
The other inclusion $s(A_1,\dots,A_m)\subseteq t(A_1,\dots,A_m)$ follows similarly from the identity $s^\ast\=\tilde{t}$.
Hence $t\=s$ holds in $\Sub\var$.
\end{proof}

\begin{cor}
Let $\var$ be a variety satisfying the generalized entropic property. Then $\Sub\var\subseteq\var$, if and only if
for every identity $t\=s$ valid in $\var$ there is a semilinear precursor for the pair $(t,s)$ valid in $\var$.
\end{cor}

\begin{cor}
Let $\var$ be an idempotent variety satisfying the generalized entropic property. Then $\var=\Sub\var$, if and only if
for every identity $t\=s$ valid in $\var$ there is a semilinear precursor for the pair $(t,s)$ valid in $\var$.
\end{cor}

\begin{cor}
Let $\var$ be an idempotent variety satisfying the generalized entropic property and assume that
$t(x_1,\dots,x_n)$ is a linear term and $s(x_1,\dots,x_n,y_1,\dots,y_m)$ is a term
such that the variables $x_1,\dots,x_n$ occur in it at most once. Then $\Sub\var$
satisfies the identity $t\=s$, if and only if $\var$ satisfies the linear identity $t\=s^\ast$.
\end{cor}

\begin{exm}
{}
\end{exm}
It follows from Theorem \ref{gl} that $\Sub\var$ satisfies all linear identities true in $\var$. This
is in accordance with Theorem \ref{gl for vs}, because
for every pair $(t,s)$ of linear terms there is a semilinear precursor $t\=s$
(indeed, $t^\ast=t$ and $s^\ast=s$), so if $t\=s$ holds in $\var$, it is satisfied in $\Sub\var$ too.

\begin{exm}
{}
\end{exm}
Let $\var$ be the variety of abelian groups. We show that $\Sub\var$ is idempotent,
i.e., $x+x\=x$ holds in $\Sub\var$, using Theorem \ref{gl for vs}.
First, we find a semilinear precursor for the pair $(x,x+x)$:
for $s(x)=x+x$ we have $s^\ast(x,y)=x+y$ and we can put $\tilde s(x)=s^\ast(x,0)$
($0$ is a constant term in any variables); indeed, $x\=x+0$ holds in $\var$.
Next, we find a semilinear precursor for the pair $(x+x,x)$: for $t(x)=x+x$ we
have $t^\ast(x,y)=x+y$, so we can substitute in $s(x)=s^\ast(x)=x$ the term $x+y$
for the variable $x$; indeed, $x+y\=x+y$ holds in $\var$.

\begin{exm}
{}
\end{exm}
Let $\var$ be the variety of entropic idempotent groupoids with $x(xy)\=y$.
We show, using Theorem \ref{gl for vs}, that $\Sub\var$ does not satisfy the identity $x(xy)\=y$.
Assume the contrary. Put $t(x,y)=x(xy)$, $s(y)=y$ and assume that
there is a semilinear precursor $t^\ast\=\tilde s$ true in $\var$. It means, there is a unary term $u$ such
that the identity $x_1(x_2y)=u(y)$ holds in $\var$. Because of idempotency,
we can assume $u(y)=y$. However, it is easy to find a groupoid in $\var$ which fails the property:
$$\begin{array}{c|ccc}
\cdot & 0& 1& 2\\ \hline
0 & 0 &1 &2\\
1 & 2 &1 &0\\
2 & 0 &1 &2
\end{array}$$

\medskip
Unfortunately, Theorem \ref{gl for vs} does not help us to decide, whether the following conjecture from \cite{P} is true.

\begin{con}\label{conj}
Let $\var$ be an idempotent variety satisfying the generalized entropic property. Then $\var=\Sub\var$, if and only if
$\var$ has a base consisting of linear identities and the identities
$f(x,\dots,x)\=x$, for all basic operations $f$.
\end{con}

Note that the backward implication is true for any idempotent variety.

All known idempotent varieties with $\var=\Sub\var$ have a linear and idempotent base.
For instance, the variety of all modes of
a given type, the variety of commutative binary modes, the variety
of differential groupoids (groupoid modes satisfying $x(yz)\=xy$),
the variety of normal semigroups (semigroup modes) and any
subvariety of this variety (in particular, varieties of
semilattices), left ($xy\=x$) and right ($xy\=y$) zero bands,
rectangular bands ($xyz\=xz$) and left ($xyz\=xzy$) and right
($zyx\=yzx$) normal bands or the variety of barycentric algebras \cite{RS3}.

We also note that at the moment we do not know any example of a non-entropic
idempotent variety $\var$ with $\var=\Sub\var$. Indeed, the only examples (known to us) of non-entropic
idempotent algebras with the generalized entropic property were shown in Example
\ref{rr}. For instance, it is straightforward to check that
$\M_9$ satisfies the identity $\underline{r}(x,\underline{r}(y,x))\=y$,
where $\underline{r}$ is any basic operation from $\underline X$, but
this identity fails in $\Sub \M_9$.

In the last example of this section we show that Conjecture \ref{conj} is false if the
assumption of idempotency is dropped.

\begin{exm}\label{var_example}
{}
\end{exm}
Consider the variety $\var$ of entropic groupoids with $(xx)y\=xy$
and $y(xx)\=yx$. Clearly, $\var$ satisfies the generalized entropic property.
It follows from Theorem \ref{gl} that $\Sub\var$ is entropic, and
it is easy to check that $\Sub\var$ satisfies the two identities.
Hence, $\Sub\var \subseteq\var$. Now, for any algebra $\A\in\var$,
we can embed $\A$ into $\Sub \A$ by $x\mapsto\{x,xx\}$ (a straightforward caculation).
Therefore, $\var=\Sub\var$. We prove that $\var$ cannot be based by linear
identities.

All identities of $\var$ are regular, i.e., have the same variables
on both sides, because the basis of $\var$ consists of regular identities.
Evidently, regular linear identities are \emph{balanced}, which means that
the number of each variable symbol, counting repetitions, is the same on both sides.
It is easy to see that consequences of balanced indentities are balanced.
Since the identities $(xx)y\=xy$ and $y(xx)\=yx$ are not balanced, they cannot be
deduced from any set of linear identities of $\var$.

\section{Stronger conjecture fails}\label{StrongHypo}

In this section we are interested in varieties that do not necessarily
possess the generalized entropic property. Our aim is to disprove an analogue of
Conjecture \ref{conj}: There is a variety generated by an idempotent
algebra $\A$ such that $\Sub \A$ exists, $\V(\Sub \A) = \V(\A)$ and $\V(\A)$ has no
base of linear and idempotent identities. The rest of the section is fully
devoted to such example.

Consider, again, the groupoid $\G_1$ from Example \ref{g1}.
$$\begin{array}{c|ccc} \cdot & a & b & c \\
\hline a & a & c & c\\
b & c & b & c \\
c & a & b & c
\end{array}$$
We already noticed that $\Sub \G_1$ exists, though $\G_1$ does not
satisfy the generalized entropic property. We show that
the groupoids $\G_1$ and $\Sub \G_1$ generate the same variety (Lemma
\ref{varieties}), but $\V(\G_1)$ has no base of linear and
idempotent identities. In fact, we prove that all linear identities
satisfied by $\G_1$ are regular (Lemma \ref{v-sets}) and thus the non-regular identities
$$(xy)x \= x\qquad\text{ and }\qquad (yx)x \= x,$$
valid in $\G_1$, are not consequences of idempotent and linear
identities of $\G_1$.

\medskip
Every term $t$ can be written in the form
$$t = t_1(t_2(\ldots t_{k-1}(t_k x)\ldots)),$$
where $t_1, \dots, t_k$ are terms and $x$ is a variable. The
variable $x$ will be called the \emph{focal of $t$} and denoted by $fc(t)$.

\begin{lm}\label{focal_nonlin}
If $\G_1$ satisfies an identity $t\=u$, then $fc(t) = fc(u)$.
\end{lm}

\begin{proof}
Assume $fc(t)\neq fc(u)$. Assign the element $c$ to $fc(t)$ and the element $a$ to all other variables of $t$ and $u$.
Then the value of $t$ is $c$ and the value of $u$ is $a$, because $a,c$ are right zeros in the subgroupoid $\{a,c\}$. Hence
$t\not\=u$ in $\G_1$.
\end{proof}

\begin{lm}\label{focal}
If $\G_1$ satisfies a linear identity $t\=u$,
$t = t_1(t_2(\ldots t_{k-1}(t_k x)\ldots))$ and
$u = u_1(u_2(\ldots u_{m-1}(u_m x)\ldots))$, then
\begin{itemize}
\item[(1)]
 $\{fc(t_i): i \leq k\}
= \{fc(u_j): j \leq m\}$. In particular, $m=k$.
\item[(2)]
For every $i \leq k$ there exists $j \leq k$ such that $t_i = u_j$ is a linear identity
of $\G_1$.
\end{itemize}
\end{lm}

\begin{proof}
To prove (1), we can assume that $fc(t) = fc(u)=x $ and there exists $y = fc(t_i)$ that does not belong
to $\{fc(u_j): j \leq m\}$. Then we assign $x = a$, $y = b$ and the rest of variables will be $c$.
It will follow that all variables in $\{fc(u_j): j \leq k\}$ will be assigned $c$, hence
all $u_j$ are equal to $c$ and $u = a$. On the other hand, $t_i = b$,
while the rest of $t_p$, $p \not = i$, are $c$.
Hence $t = c$ and $t \not = u$ under such assignment of variables.

To show (2), for any $t_i$ we pick $u_j$ with the same focal $y$.
Suppose that $t_i \not = u_j$ for some assignment of variables. Then $y$ is assigned to
$a$ or $b$.

If $y=a$, then $\{t_i,u_j\} = \{a,c\}$ under such assignment. Say,
$t_i = a$ and $u_j = c$. Let $fc(t)=fc(u)$ be assigned to $b$ and
all $fc(t_p)$, $p\neq i$, and $fc(u_q)$, $q\neq j$, to $c$. Under
such assignment we get that $t=c$ and $u = b$, a contradiction
with $t=u$ in $\G_1$.

The case of $y=b$ is shown similarly by interchanging $a$ and $b$.
\end{proof}

\begin{lm}\label{v-sets}
Every linear identity of $\G_1$ is regular.
\end{lm}

\begin{proof}
Let $r(t,u)$ be the number of distinct variables in the identity $t\=u$
(e.g., $r(xy,(xz)y)=3$).
We argue by induction on $r(t,u)$. If $r(t,u) = 1$, then $t=u=x$ and the statement is true.

Suppose we know that every linear identity $t'\=u'$ with $r(t',u') \leq n$ is regular and
consider a linear identity $t\=u$ with $r(t,u)=n+1$. Then, according to Lemma \ref{focal},
$t = t_1(t_2(\ldots t_{k-1}(t_k x)\ldots))$ and
$u = u_1(u_2(\ldots u_{k-1}(u_k x)\ldots))$, for some $k$
and some terms $t_i,u_i$ such that for every $i \leq k$ there exists $j \leq k$ with $t_i \= u_j$
satisfied in $\G_1$. This is indeed a linear identity and $r(t_i,u_j) \leq n$, because $x$ does not
occur in $t_i\=u_j$. By induction hypothesis, $t_i\=u_j$ is regular. Hence the set of variables occuring
in $t$ is a subset of the set of variables occuring in $u$. Similarly, applying Lemma \ref{focal} on the
identity $u\=t$, we obtain that the latter set is a subset of the former one. Consequently, the identity
$t\=u$ is regular.
\end{proof}

As a byproduct we also get a description of linear identities satisfied in $\G_1$.
For this, we define \emph{focally equivalent} terms $t\equiv_fu$, recursively by the length of $t,u$:

\begin{itemize}
\item[(1)] If one of $t,u$ has only one variable $x$ then $t\equiv_fu$ if and only if $t=u=x$.
\item[(2)] If both $t,u$ have more than one variable and $t = t_1(t_2(\ldots t_{k-1}(t_k x)\ldots))$,\\
$u= u_1(u_2(\ldots u_{l-1}(u_l y)\ldots))$, then $t\equiv_fu$ if and only if $k=l$, $x = y$,
for every $i\leq k$ there is $j\leq k$ such that $t_i\equiv_f u_j$
and for every $i\leq k$ there is $j\leq k$ such that $u_i\equiv_f t_j$.
\end{itemize}

\begin{cor}\label{linearG1}
$\G_1$ satisfies a linear identity $t\=u$, iff the terms $t,u$ are focally equivalent.
\end{cor}

\begin{proof}
Apply induction and Lemmas \ref{focal_nonlin} and \ref{focal}.
\end{proof}

\begin{lm}\label{varieties}
$\G_1$ and $\Sub \G_1$ generate the same variety.
\end{lm}

\begin{proof}
Since an idempotent algebra always embeds into its algebra of subalgebras (provided it exists),
it is sufficient to find an embedding of the groupoid $\Sub \G_1$ into the product
$\G_1 \times \G_1$.
We notice that there are two homomorphisms from $\Sub \G_1$ onto $\G_1$:

\begin{align*}
f_1(\{a\}) = f_1(\{a,c\})=f_1(\{a,b,c\}) &= a,\\
f_1(\{b\}) &= b,\\
f_1(\{c\}) =f_1(\{b,c\}) &= c
\end{align*}
and
\begin{align*}
f_2(\{b\}) = f_2(\{b,c\})=f_2(\{a,b,c\}) &= b,\\
f_2(\{a\}) &= a,\\
f_2(\{c\})=f_2(\{a,c\}) &= c.
\end{align*}
It is easy to check that $ker(f_1) \cap ker(f_2) = 0$, hence $\Sub \G_1$
is a subdirect power of $\G_1$.
\end{proof}

{\bf Acknowledgments.} The paper was initiated during the first
author's visit to Warsaw University of Technology in summer of
2003. The support and welcoming atmosphere created by Prof.
A.Romanowska and the group of her collaborators are greatly
appreciated. The work on this project was further inspired by
the INTAS workshop at Charles University in Prague in summer of 2004
and the visits of the last and first authors to Warsaw in 2005.

Several times we helped us with the automated theorem prover Otter
\cite{otter}. A number of proofs in Section 3 are translated and
simplified versions of Otter's ones.

We want to thank Peter Jones for his interest in this work, in
particular providing us with an alternative direct proof of
Proposition \ref{bands}. Many helpful suggestions were also
communicated to us by Jonathan Smith, Michal Stronkowski and
George Gr\"atzer.

\end{document}